\newcommand{\ot}{\otimes}

\newcommand{\ra}{\rightarrow}

\newcommand{\BR}{\mathbb{R}}

\newcommand{\BZ}{\mathbb{Z}}

\newcommand{\Id}{\mathbf{Id}}

\newcommand{\Spinc}{\ensuremath{\mathrm{spin}^c}}

\documentclass[twoside,12pt]{amsart}

\usepackage{latexsym}

\usepackage{amsmath}

\usepackage{amsthm}

\usepackage{amscd}

\usepackage{xypic} 

\usepackage{amsfonts}

\usepackage{amssymb}

\usepackage{graphicx}

\newtheorem{df}{Definition}

\newtheorem{thm}[df]{Theorem}

\newtheorem{cor}[df]{Corollary}

\newtheorem{lem}[df]{Lemma}

\newtheorem{prop}[df]{Proposition}

\newtheorem{quest}{Question}

\begin{document}

\title[4-Manifolds with free circle actions]{Seiberg-Witten invariants of
4-manifolds with free circle actions}

\author{Scott Baldridge}

\address{Department of Mathematics, Michigan State University \newline
\hspace*{.375in}East Lansing, Michigan 48824}

\email{\rm{baldrid1@pilot.msu.edu}}

\date{September, 1999}

\maketitle




\section{Introduction}

The main result of this paper describes a formula for the Seiberg-Witten
invariant of a $4$-manifold $X$ which
admits a nontrivial free $S^1$-action.
A free circle action on $X$ is classified by its orbit space, a  3-manifold
$M$, and its Euler class $\chi \in
H^2(M;\BZ)$. If $\chi=0$, then $X=M \times S^1$, and it is well-known that
the Seiberg-Witten
invariants of $X$ are equal to the $3$-dimensional Seiberg-Witten
invariants of $M$.

 Our result expresses the Seiberg-Witten invariants of $X$ are in terms of
the Seiberg-Witten invariants of
$M$ and the Euler class $\chi$:

\begin{thm}
Let $X$ be a smooth $4$-manifold with $b_+ \geq 2$ and a free circle
action.  Let $M^3$ be the smooth orbit
space and suppose that the Euler class $\chi \in H^2(M;\BZ)$ of the free
circle action is not torsion.  Let
$\xi$ be a \Spinc\ structure over $X$.  If $\xi$ is not pulled up via
$\pi:X\ra M$, then $SW_X(\xi) = 0$.
Otherwise, let $\xi^*$ be a \Spinc\ structure on $M$ such that $\xi =
\pi^*(\xi^*)$, then
\begin{eqnarray}
SW_{X}^4(\xi) & = & \sum_{\xi' \equiv \xi^*\mod{\chi}}  SW_{M}^3(\xi').
\end{eqnarray}
\label{thm:maintheorem}
\end{thm}

The difference of two \Spinc\ structures gives rise to a well-defined
element $\xi'-\xi \in H^2(X;\BZ)$. For more
information, see section~(\ref{sec:spinc}). Because $\chi$ is nontorsion, the
equivalence relation in the above theorem is well-defined.  The pullback of
a \Spinc\ structure is discussed in
section~(\ref{sec:pullback}).

As an application of this theorem we shall produce a nonsymplectic
$4$-manifold with a free circle action whose
orbit space fibers over $S^1$. This example runs counter to intuition since
there is a well-known conjecture of Taubes
that
$M^3\times S^1$ admits a symplectic structure if and only if $M^3$ fibers
over the $S^1$.
Furthermore, there is
evidence~\cite{symp:sympmanfreecircleaction}  which suggests that many such
$4$-manifolds are, in fact, symplectic.
As another application of our formula, we construct a $3$-manifold which is
not the orbit space of any symplectic
$4$-manifold with a free circle action.
A corollary of the main theorem is a formula for the Seiberg-Witten invariant of the total space of a circle 
bundle over a surface.  This formula can be thought of as the $3$ dimensional analog of the $4$ dimensional 
formula.




\section{Classifying free circle actions}

Let $X$ be an oriented connected $4$-manifold carrying a smooth free
$S^1$-action.  Its orbit space $M$ is a
$3$-manifold whose orientation is determined, so that, followed by the
natural orientation on the orbits, the
orientation of $X$ is obtained. Choose a smooth connected loop $l$
representing the the Poincar\'e dual
$PD(\chi)\in H_1(M;\BZ)$.   Remove a tubular neighborhood $N \cong D^2
\times l$ of $l$
from $M$, and  set $X_0 = (M \setminus N)\times S^1$. View $X_0$ as an
$S^1$-manifold whose action is given by rotation in
the last factor. Let
$m$ be the meridian of
$l$, and let
$t$ be an orbit in
$X_0$.  We then have:
\begin{lem}
The manifold $X$ is diffeomorphic (by a bundle isomorphism) to the manifold
\begin{eqnarray}
X(l) = X_0 \cup_\varphi D^2 \times T^2
\end{eqnarray}
where $\varphi:T^3 \ra \partial X_0$ is an equivariant diffeomorphism which
evaluates $\varphi_*([\partial(D^2
\times pt)] = [m+t]$ in homology.
\end{lem}
When gluing $D^2 \times T^2$ into the boundary of a manifold, the resulting
closed manifold is determined up to
diffeomorphism by the image in homology of $[\partial (D^2 \times pt)]$.
(For example,
see~\cite{sw:productformulaalongT3}.)
\begin{proof}
The manifold $X$ is a principal $S^1$-bundle. Since $\chi$ evaluates on any
2-cycle in $M \setminus N$ by intersecting
that $2$-cycle against $l$, it follows that
the restriction of the Euler class $\chi$ restricts trivially to $M
\setminus N$.
Therefore, the $S^1$-bundle is trivial over $M \setminus N$, and
$\pi^{-1}(M \setminus N)$ is diffeomorphic to
$X_0$.  Similarly, $\pi^{-1}(N)$ is diffeomorphic to $D^2 \times
S^1 \times S^1$.  Let $m'$, $l'$,
and $t'$ be the circles which correspond to the factors in $D^2 \times S^1
\times S^1$ respectively.

Construct a manifold $X(l)$ as above using a bundle isomorphism
$\varphi:\partial (D^2 \times S^1) \times S^1
\ra X_0$.  Bundle isomorphisms covering the identity are classified up to
vertical equivariant isotopy by
homotopy classes of maps in $[\partial(D^2 \times S^1), S^1 ] =  \BZ \oplus
\BZ$.  Explicitly, an equivariant
map $\varphi$ inducing $1_{\partial(D^2 \times S^1)}$ is classified by
integers (r,s) where $\varphi_*[m']=[m] +
r[t]$ and $\varphi_*[l']=[l]+s[t]$.   A bundle automorphism $\Phi$ of $(D^2
\times S^1) \times S^1$ can be
constructed such that $\Phi_*[m']=[m']$ and $\Phi_*[l']= [l']+s[t']$ for any
$s \in \BZ$.  These bundle
automorphisms are just the equivariant maps classified by $[D^2 \times S^1,
S^1]=H^1(D^2 \times S^1;\BZ)$.
Therefore the resulting bundle $X(l)$ depends only on the integer $r$ and
the homology class $[l]$.  In
particular, the obstruction to extending the constant section
\[ M\setminus N\to X_0=(M\setminus N)\times S^1 \]
over $D^2 \times S^1$ lies in $H^2(D^2 \times S^1, \partial (D^2 \times
S^1); \BZ)$ and is given by $r$.  The Euler class of $X(l)$ is then
$PD(r[l]) = r\chi$.  Taking $r=1$ produces
the desired bundle.
\end{proof}

 From now on we shall work with $X(l)$ and refer to it as $X$.  Furthermore,
it is clear from
the construction above that the map $\varphi$ can be chosen so that in
homology,
\begin{eqnarray}
\varphi_* = \left(\begin{array}{ccc} 1 & 0 & 0 \\ 0 & 1 & 0\\ 1 & 0 & 1
\end{array}\right)
\end{eqnarray}
with respect to the basis $\{ [m], [l], [t]\}$.
\label{classifycircleact}




\section{Gluing along $T^3$}

Since we have  $X=X_0\cup_{\varphi}(D^2\times T^2)$ we may apply the gluing
theorem of Morgan, Mrowka, and
Szab\'{o} ~\cite{sw:productformulaalongT3}. Recall that
$\varphi_*([m'])=[m+t]$.
\begin{thm}[Morgan, Mrowka, and Szab\'{o}]
If the \Spinc\ structure $\xi$ over $X$ restricts nontrivially to
$D^2\times T^2$, then $SW_X(\xi)=0$. For each
\Spinc\ structure $\xi_0 \ra X_0$ that restricts trivially to $\partial
X_0$, let $V_{X}(\xi_0)$ denote the set
of isomorphism classes of \Spinc\ structures over $X$ whose restriction to
$X_0$ is equal to $\xi_0$.  Then we
have
\begin{equation}
\small
\sum_{\xi\in V_{X}(\xi_0)} SW_{X}(\xi)
 =
\sum_{\xi\in V_{M\times S^1}(\xi_0)} SW_{M\times S^1}(\xi) +
\sum_{\xi\in V_{X_{0/1}}(\xi_0)} SW_{X_{0/1}}(\xi),
\label{eq:main_eq}
\end{equation}
where the manifold $X_{0/1} = X_0 \cup_{\varphi_{0,1}} D^2\times T^2$ is
defined by the map $\varphi_{0,1}$
which maps $[m'] \mapsto [t]$ in homology.
\end{thm}

In our situation, this formula simplifies significantly. Let $i$ denote
the inclusion of $\partial X_0$ into $X_0$.  A study of the long exact
sequences in homology shows that the left
hand side consists of a single term when $i_*[m+t]$ is indivisible.  Since
$i_*[t]$ is independent of $i_*[m]$
and $i_*[t]$ is a primitive class in $H_1(X_0;\BZ)$, $i_*[m+t]$ is such a
class.  Therefore, the formula enables
the calculation of the SW invariants of $X$ in terms of the SW invariants
of $M\times S^1$ and  a manifold
$X_{0/1}$.

The manifold $X_{0/1}$ admits a semi-free $S^1$-action whose fixed point
set is a torus. Its orbit space
is $M\setminus N$, and $\partial(M\setminus N)=\partial N$ is the image of
the fixed point set.  The condition
$b_+(X)\geq 2$ of the main theorem implies that
$b_+(X_{0/1})>1$ and that
\[ \mathrm{rank}\; H_1(M
\setminus N, \partial (M \setminus N);\BZ) > 1.  \]
The two statements are proved as follows.  The Gysin sequence 
\begin{equation}
\xymatrix{
H^2(M;\BZ) \ar[r]^{\pi^*} & H^2(X;\BZ) \ar[r] & H^1(M;\BZ) \ar[r]^{\cdot \cup \chi} & 
H^3(M;\BZ)}
\label{eq:gysin_seq}
\end{equation}
implies 
\begin{multline}
\label{eq:secondcohomologyofX}
H^2(X;\BZ) \cong \left( H^2(M;\BZ) / <\chi>\right) \\ \oplus \ker \left(\cup\chi :
 H^1(M;\BZ) \ra H^3(M;\BZ)\right).
\end{multline}
Each component of the direct sum above has rank $b_1(M) -1$.  The bilinear form of 
$X$ is the direct sum of hyperbolic pairs which implies that $b_+(X) = b_1(M) -1$. Since $[l]$ is not a torsion 
element, removing $N$ 
from $M$ implies the rank of $H_1(M\setminus N, \partial(M\setminus N);\BZ)$ is also 
$b_1(M)-1$. The second statement now follows because $b_1(M) - 1 = b_+(X) >1$. The first statement requires the 
following Mayer-Vietoris sequence
\[
H_3(T^3;\BZ) \ra  H_2(X_0;\BZ) \oplus 
H_2(D^2\times T^2;\BZ) \ra  H_2(X_{0/1};\BZ) \stackrel{0}{\ra}  H_1(T^3;\BZ). \]
The rank of $H_2(X_0;\BZ)$ is $2b_1(M)-1$ and the rank of the image of the first map is $2$. Therefore 
$b_2(X_{0/1}) = 2b_1(M) - 2$.  Since the bilinear form of $X_{0/1}$ is also a direct sum of hyperbolic pairs, 
$b_+(X_{0/1}) > 1$.

\begin{prop}\label{vanishing}
Let $X$ be a smooth closed oriented $4$-manifold with a smooth semi-free
circle action and
$b_+(X)>1$.  Let $X^* = X/S^1$ be its orbit space.  Suppose that $X^*$ has
a nonempty boundary and
$\mathrm{rank} \; H_1(X^*,\partial X^*;\BZ)>1$.  Then $\mbox{SW}_X \equiv 0$.
\end{prop}
\begin{proof}
Let $F$ denote the fixed point set of $X$ and $F^*$ its image in $X^*$.
Then $\partial X^*\subset F^*$.
The restriction of the circle action to $X\setminus F$ defines a principal
$S^1$-bundle whose Euler class lies in
$H^2(X^*\setminus F^*;\BZ)$. Let $\chi'\in H_1(X^*,F^*;\BZ)$ denote its
Poincar\'e dual. Consider the exact sequence
\begin{multline*} 0\to H_1(X^*,\partial X^*;\BZ) \xrightarrow{i_*}
H_1(X^*,F^*;\BZ)\to \\
\to H_0(F^*,\partial X^*;\BZ)\to H_0(X^*,\partial X^*;\BZ).
\end{multline*}
Since the rank of $ H_1(X^*,\partial X^*;\BZ)$ is greater than 1, there is
a class in $i_*(H_1(X^*,\partial X^*;\BZ))$
which is primitive and not a multiple of $\chi'$.
This class may be represented by a path $\alpha$ in $X^*$ which
starts and ends on $\partial X$ but is otherwise disjoint from $F^*$.

The preimage $S=\pi^{-1}(\alpha)$ is a 2-sphere of self-intersection $0$ in
$X$. The Gysin sequence gives:
\[ H_3(X^*,F^*,\BZ)\to H_1(X^*,F^*,\BZ)\xrightarrow{\rho} H_2(X,F,\BZ)\to
H_2(X^*,F^*,\BZ) \]
where $\rho_*(i_*[\alpha])=[S]$. The image of $H_3(X^*,F^*,\BZ)\cong\BZ$ in
$H_1(X^*,F^*,\BZ)$ is generated by $\chi'$.
Since $i_*[\alpha]$ is primitive and not a multiple of $\chi'$, the class
$[S]\in {\text{Im}}\rho\subset H_2(X,F,\BZ)$
is not torsion; hence $[S]$ is nontorsion as an element of $H_2(X;\BZ)$.

It now follows from \cite{sw:class_turkey} that $\mbox{SW}_X \equiv 0$.
\end{proof}

This type of vanishing theorem is quite common for 4-manifolds with circle
actions.  For instance, it follows from
\cite{circact:Circ_act_on_four_man} that Seiberg-Witten invariants vanish
for simply connected 4-manifolds
which have $b_+ >1$ and a smooth circle action.

Proposition~\ref{vanishing} implies that the
 formula~(\ref{eq:main_eq}) simplifies to
\begin{eqnarray}
SW_{X}(\xi) \! &=&  \sum_{\xi'\in V_{M \times S^1}(\xi|_{X_0})} \;SW_{M
\times S^1}(\xi').
\label{eq:equation_with_MT}
\end{eqnarray}




\section{Understanding the $\mathrm{spin}^c$ structures}

In this section we shall prove that all basic classes of $X$ come from
\Spinc\ structures that are pulled up
from $M$ (in a suitable sense).  We shall also identify the \Spinc\
structures in the set $V_{M\times
S^1}(\xi|_{X_0})$ coming from the gluing theorem.

\subsection{$\mbox{Spin}^c$ structures}
\label{sec:spinc}

First recall some basic facts about \Spinc\ structures. The set of \Spinc\
structures lifting the frame bundle
of a $4$-manifold $X$ is a principal homogeneous space over
$H^2(X;\BZ)$: given two \Spinc\
structures $\xi_1,\xi_2$ their difference $\delta(\xi_1,\xi_2)$ is a
well-defined element of $H^2(X;\BZ)$. For
details, see ~\cite{sw:alg_surface} or~\cite{gauge:spingeometry}.

Likewise, if $\xi$ is a \Spinc\ structure and $e\in H^2(X;\BZ)$ is a
2-dimensional cohomology class, there is a
new \Spinc\ structure $\xi+e$.  Let $W_\xi$ be spinor bundle associated
with $\xi$, then the new spinor bundle
is $W_\xi \ot L_e$ where $L_e$ is the unique line bundle with first Chern
class $e$.

For all \Spinc\ structures, a line bundle $L_\xi$ can be associated to
$\xi$ called the determinant line bundle.
Let $(\xi,L_\xi)$ be a pair consisting a \Spinc\ structure $\xi$ whose
determinant line bundle is $L_\xi$.
Given two \Spinc\ structures $(\xi_1,L_1), \; (\xi_2,L_2)$, the difference
of their determinant line bundles is
$c_1(L_1)-c_1(L_2) = 2 e$ for some element $e \in H^2(X;\BZ)$.  If
$H^2(X;\BZ)$ has no $2$-torsion, then $e$ is
well-defined and $c_1(L_\xi)$ determines the \Spinc\ structure for $(\xi,
L_\xi)$.  When $H^2(X;\BZ)$ has
$2$-torsion, one has a choice of two or more possible square roots of $2 e$
and it seems that $e$ is not
well-defined.  However, the difference element $\delta(\xi_1,\xi_2)$
satisfies $c_1(L_1)-c_1(L_2) = 2
\delta(\xi_1,\xi_2)$ and so there is a unique element in $H^2(X;\BZ)$ which
determines the difference of two
\Spinc\ structures even in the presence of $2$-torsion.  So while
$c_1(L_\xi)$ does not determine $\xi$ in this
case, the difference between two \Spinc\ structures is still well-defined.

\subsection{Pullbacks of $\mbox{spin}^c$ structures}
\label{sec:pullback}

The \Spinc\ structures on a 3-manifold $M$ are defined by a
pair $\xi = (W, \rho)$ consisting of a rank 2 complex bundle $W$
with a hermitian metric (the spinor bundle) and an action $\rho$ of
1-forms on spinors,
\[\rho: T^*M \ra \mbox{End}(W),\]
which satisfies the following property
\[\rho(v)\rho(w) + \rho(w)\rho(v) = -2<v,w>\Id_W.\]

For a 4-manifold the definition is similar, but  consists
of a rank 4 complex bundle with an action on the cotangent space
that satisfies the same property. There is a natural way to define the
pullback of a \Spinc\ structure.
Let $\eta$ denote the connection 1-form of the circle bundle $\pi: X\to M$,
and let $g_M$ be a metric on  $M$, then we
can endow
$X$ with the metric $g_{X} =  \eta \ot \eta + \pi^*(g_M)$.  Using this
metric, there is an orthogonal splitting
\[T^*X \cong  \BR\eta \oplus \pi^*(T^*M).\]
If $\xi=(W,\rho)$ is a \Spinc\ structure over $M$, define the pullback of
$\xi$ to be $\pi^*(\xi) = (\pi^*(W)
\oplus \pi^*(W), \sigma)$ where the action
\[ \sigma: T^*X \ra \mbox{End}(\pi^*(W) \oplus \pi^*(W))\]
is given by
\[ \sigma (b\eta + \pi^*(a)) = \left( \begin{array}{cc} 0 & \pi^*(\rho(a))
+ b\Id_{\pi^*}(W) \\
\pi^*(\rho(a)) - b\Id_{\pi^*(W)} & 0 \end{array} \right). \]
One can easily check that this defines a \Spinc\ structure on $X$. Note
that the
first Chern class of $\pi^*(\xi)$ is just $\pi^*(c_1(L_{\xi}))$.
The other pulled back \Spinc\
structures are now obtained by the addition of classes $\pi^*(e)$ for $e\in
H^2(M;\BZ)$.

There are \Spinc\ structures on $X$ which do not arise from \Spinc\
structures that are pulled up from $M$.
In the next section we show that the Seiberg-Witten invariants vanish for
these \Spinc\ structures.

\subsection{$\mbox{Spin}^c$ structures which are not pullbacks}

Fix a \Spinc\ structure $\xi_0 = (W_0,\rho)$ on $M$ and consider its
pullback $\xi =\pi^*(\xi_0)$ over $X$.  
Looking at the Gysin sequence~(\ref{eq:gysin_seq}), if a class 
$e\in H^2(X;\BZ)$ is not in the image of $\pi^*$, then $\xi
+ e$ is not a \Spinc\ structure which is pulled back from $M$.

\begin{lem}
If $(\xi, L_\xi)$ is a \Spinc\ structure on $X$ which is not
pulled back from $M$,  then $SW_{X}(\xi) = 0 $.
\end{lem}
\begin{proof}
We claim that there exists an embedded torus which pairs nontrivially with
$c_1(L_\xi)$.  Then by the
adjunction inequality~\cite{sw:adjunctioninequality} the \Spinc\ structure
$\xi$ has Seiberg-Witten invariant
equal to zero.  Let 
\[\mathbf{H}=\ker (\cdot \cup \chi : H^1(M;\BZ) \ra
H^3(M;\BZ))\]
in equation~(\ref{eq:secondcohomologyofX}), and consider for a moment
the projection of $c_1(L_\xi)$ onto the first factor of $\mathbf{H}\oplus
\pi^*(H^2(M;\BZ))$ by changing the
\Spinc\ structure by an element of $\pi^*(H^2(M;\BZ))$. Since $\xi$ is not
pulled back from $M$,
$c_1(L_\xi)|_\mathbf{H}\ne 0$, and since
$H^1(M;\BZ)$ is a free abelian group, $c_1(L_\xi)|_\mathbf{H}$ is not a
torsion class.

Examining the Gysin sequence,  $c_1(L_\xi)|_\mathbf{H} \in H^2(X;\BZ)$ maps
to a class $\beta \in H^1(M;\BZ)$,
$\beta\cup\chi=0$. Thus the Poincar\'e
dual of $\beta$ can be represented by a surface $b$, and there is a
$1$-cycle $\lambda$ in
$M\setminus N\  {\text{rel}}\ \partial$ such that $[\lambda] \cdot [b]
\not= 0$. Since
$\partial N$ is connected, $[\lambda]$ is actually represented by a loop
$\lambda$ in $M\setminus N$.  The
preimage $\pi^{-1}(\lambda) = \lambda \times S^1$ in $X$ is a torus, and
$c_1(L_\xi)|_\mathbf{H}\cdot
[\pi^{-1}(\lambda)] = [b]\cdot[\lambda]\ne 0$.

On the other hand, if $A\in \pi^*H^2(M;\BZ)$ then it Poincar\'e dual is
represented by a loop $\alpha$ in $M$ which may
be chosen disjoint from $\lambda$. Thus $A\cdot [\pi^{-1}(\lambda)] =
0$. This means that
$c_1(L_\xi)\cdot [\pi^{-1}(\lambda)]\ne 0$, as required.
\end{proof}

\subsection{Identifying the set $V_{M \times S^1}(\xi|_{X_0})$}

According to the previous lemma, the only nontrivial Seiberg-Witten \Spinc\
structures are those pulled up
from  $M$. Thus far we have seen that for such a \Spinc\ structure
$\xi=\pi^*(\xi^*)$ with $\xi_0= \xi|_{X_0}$, we have
\[ SW_{X}(\xi) = \sum_{\xi'\in V_{M\times S^1}(\xi_0)} SW_{M\times S^1}(\xi').\]
Let $\tilde{\pi}:M \times S^1 \ra M$ be the projection. We
identify the set $V_{M\times S^1}(\xi_0)$ of isomorphism classes of \Spinc\
structures over $M\times S^1$ which restrict
on $X_0$ to $\xi_0$.

\begin{lem}
$V_{M \times S^1}(\xi_0)= \{ \; \tilde{\pi}^* \left(\xi^* + n \cdot
\chi\right) \;|\; n \in \BZ  \;\}$.
\end{lem}
\begin{proof}
The diagram
\[ \xymatrix{
X \ar[rddd]_\pi & X_0 \ar[l]_{inc} \ar[r]^{inc}
\ar[d]^{\tilde{\pi}|_{M\setminus N}} & M\times S^1 \ar[lddd]^{\tilde{\pi}} \\
& M\setminus N \ar[dd]^{inc} \\
\\
& M
}
\]
induces \Spinc\ structures on $X$, $X_0$, and $M \times S^1$ which satisfy
\[inc^*(\pi^*(\xi^*)) = \xi_0 = inc^*(\tilde{\pi}^*(\xi^*).\]
Recall that $\xi$ is the only \Spinc\ structure induced on $X$ by $\xi_0$
since $i_*[m+t]$ is indivisible.
Since  $\tilde{\pi}^*(\xi^*) \in V_{M \times S^1}(\xi_0)$, the set of
\Spinc\ structures on $M \times S^1$
is $\{\tilde{\pi}^*(\xi^*) + e| e \in H^2(M \times S^1; \BZ)  \}$.  Now
$\tilde{\pi}^*(\xi^*)+e$ lies in
$V_{M \times S^1}(\xi_0)$ if and only if $inc^*(\pi^*(\xi^*)+e) =\xi_0$,
i.e. if and only if $inc^*(e)=0$.
Therefore,
\begin{eqnarray}
V_{M \times S^1}(\xi_0) = \left\{ \tilde{\pi}^*(\xi^*) + e \; |\;  inc^*(e)
= 0 \right\}.
\label{eq_VMXS}
\end{eqnarray}

The kernel of $inc^*$ is equal to the image of $j^*$ in the diagram below.
\[ \xymatrix{
H^2(M \times S^1, (M \setminus N)\times S^1; \BZ) \ar[r]^{j^*} \ar[d]^{PD} & H^2(M \times
S^1;\BZ) \ar[r]^{inc^*} \ar[d]^{PD} & H^2(X_0;\BZ) \ar[d]^{PD} \\
H_2(D^2 \times T^2;\BZ) \ar[r] & H_2(M \times S^1;\BZ)
\ar[r] & H_2(X_0,\partial X_0;\BZ)\\
n[T^2] \ar[r]_{j_*} & n[l \times t] \ar[r] & 0 \\
}\]
However $j_*[{\text{pt}}\times T^2] =[l\times t]$, and since
$\tilde{\pi}^*(\chi) = PD^{-1}[l \times t]$, the lemma
follows.
\end{proof}

\subsection{Relationship between $SW^3$ and $SW^4$}
\label{sec:donald}

The following is a well-known fact about the relationship between the
$3$-dimensional Seiberg-Witten invariants
and the $4$-dimensional invariants.

\begin{prop}[cf. Donaldson~\cite{sw:swand4man}]
After making a suitable choice of orientations for $M$ and $M\times S^1$,
the following equality holds
\[SW^3_M(\xi) = SW^4_{M\times S^1}(\tilde{\pi}^*(\xi))\]
for a \Spinc\ structure $\xi$ over $M$.
\label{prop:donaldson}
\end{prop}

A natural choice of orientations for $M \times S^1$ and $M$ is induced by
the orientation of the circle action
on $X$. This completes the proof of Theorem~\ref{thm:maintheorem}.




\section{Applications and examples}

\subsection{An application}

An immediate corollary to the main theorem is the calculation of the $3$ dimensional Seiberg-Witten 
invariants for the total space of a circle bundle over a surface. The following corollary can also be
derived from~\cite{sw:sw_inv_seifert_space} using different techniques.

\begin{cor}
Let $\pi:Y\ra \Sigma_g$ be a smooth $3$-manifold which is the total space of a circle bundle over a surface of 
genus 
$g>0$.  Let $c_1(Y) = n\lambda \in H^2(\Sigma_g;\BZ)$ where $\lambda$ is the generator.  The only invariants 
which are not zero on $Y$ 
come from \Spinc\ structures which are pulled back $\pi:Y \ra \Sigma_g$.  Hence,
\[ SW_Y(\pi^*(s\lambda)) = \sum_{t \equiv s\mod{n} } SW_{\Sigma_g \times S^1} (\tilde{\pi}^*(t\lambda)) \]
where $\tilde{\pi}:\Sigma_g \times S^1 \ra \Sigma_g$.
\end{cor}

\begin{proof}
Let $\pi:Y \ra \Sigma_g$ be the total space of a circle bundle over $\Sigma$ with Euler class $n\lambda$.
Then the manifold $Y \times S^1$ can be thought of as a smooth $4$-manifold with a free circle action
which orbit space is $\Sigma_g \times S^1$.  The Euler class of the action is 
$\tilde{\pi}^*(n\lambda))$.  Applying the main theorem gives
\[SW_{Y \times S^1}^4((\pi,id)^*(\tilde{\pi}^*(s\lambda))) = \sum_{\tilde{\pi}^*(t\lambda) \equiv 
\tilde{\pi}^*(s\lambda)\mod{\tilde{\pi}^*(n\lambda )}}  SW_{\Sigma \times S^1}^3(\xi')\]
the right hand side of the equation.  Applying Proposition~\ref{prop:donaldson} shows that $SW^4=SW^3$ in this 
case.
\end{proof}

Combining the Seiberg-Witten polynomial for the product of a surface with a circle, 
\[\mathcal{SW}_{\Sigma_g \times S^1}(t) = (t^1 - t^{-1})^{2g-2},\]
with the previous results gives a formula for the Seiberg-Witten polynomial in terms of the Euler
class and the genus of the surface.

\begin{cor}
Let $\pi:Y \ra \Sigma_g$ be the total space of a circle bundle over a genus $g$ surface.
Assume $c_1(Y) = n \lambda$ where $\lambda \in H^2(\Sigma_g;\BZ)$ is the generator and
$n$ is an even number $n=2l \not= 0$, then the Seiberg-Witten 
polynomial of $Y$ is
\begin{eqnarray*}
\mathcal{SW}_Y(t) & = & \mbox{sign}(n) \sum_{i=0}^{|l|-1} \sum_{k=-(2g-2)}^{k=2g-2} (-1)^{(g-1) +i +k|l|} 
\binom{2g-2}{(g-1) +i +k|l|}  t^{2i}
\end{eqnarray*}
where $t=\mbox{exp}(\pi^*(\lambda))$ and defining the binomial cofficient $\binom{p}{q} = 0$ for $q<0$ and 
$q>p$.  For the formula where $n$ is odd, replace $l$ by $n$ and $t^{2i}$ by $t^i$.
\end{cor}

If one uses \cite{sw:swequalmilnortorsion} to calculate the Milnor torsion for a circle bundle $Y$ over a 
surface, one finds that the invariant is identically $0$.  This is because all \Spinc\ 
structures on $Y$ with nontrivial 
invariants have torsion first Chern class. Turaev introduced another type of torsion in 
\cite{sw:tor_inv_spinc_on_three, sw:comb_form_for_sw_of_three} and a 
combinatorially defined function on the set of \Spinc\ structures $T:\mathcal{S}(Y) \ra \BZ$ derived 
from this torsion, and showed that this function 
was the Seiberg-Witten 
polynomial up to sign. Hence, principal $S^1$-bundles over surfaces provide simple examples which illustrate the 
difference between Milnor torsion and Turaev torsion.

\subsection{A construction and a calculation}

The following construction is similar to but simpler than  the main
construction
in~\cite{sw:knots_links_and_four_man}.  Let $Y_K$ denote the manifold
resulting from  $0$-surgery on a knot $K$
in $S^3$.  Let $m$ be a meridian of the knot in $Y_K$.  Let $m_1, m_2, m_3$
be loops that correspond to the
$S^1$ factors of $T^3$.  Construct a new manifold
\[ M_K = T^3 \#_{m_1 = m}Y_K = [T^3 \setminus (m_1 \times D^2)] \cup [Y_K
\setminus (m \times D^2)] \]
by removing tubular neighborhoods of $m$ and $m_1$ and fiber summing the
two manifolds along the boundary such
that $m= m_1$ and such that $\partial D^2$ is sent to $\partial D^2$.

This is a familiar construction. If one forms a link $L$ from the Borromean
link by taking the composite of the
first component with the knot $K$ (see Figure 1), then $M_K$ is the result
of surgery on $L$ with each surgery
coefficient equal to $0$. If $K$ is a fibered knot, then the resulting
manifold $T^3\#_{m_1 = m} Y_K$ is a
fibered 3-manifold.

\begin{figure}
\includegraphics{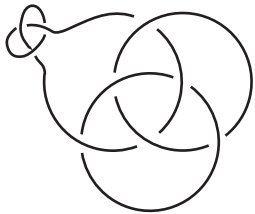}
	\caption{$M_K$ before surgery. }
	\label{fig:borromean1}
\end{figure}

Consider the formal variables $t_{\beta} = \mbox{exp}({PD(\beta)})$ for
each $\beta \in H_1(M;\BZ)$ which
satisfy the relation $t_{\alpha + \beta} = t_{\alpha}t_{\beta}$.  The
Seiberg-Witten polynomial $\mathcal{SW}$
of $X$ is a Laurent polynomial with variables $t_{\beta}$ and coefficients
equal to the Seiberg-Witten invariant
of the \Spinc\ structure defined by $t_\beta$.

\begin{thm}[Meng and Taubes~\cite{sw:swequalmilnortorsion}]
In the situation above
\begin{eqnarray}
\mathcal{SW}^3_{M_K} = \Delta_K(t_{m_1}^2)
\end{eqnarray}
where $\Delta_K$ is the symmetrized Alexander polynomial of $K$.
\label{thm:mengtaubes}
\end{thm}

For example, the manifold $M_K$ in Figure~\ref{fig:borromean1} where $K$ is
the trefoil knot has Seiberg-Witten
polynomial
\[ \mathcal{SW}^3_{M_K}(t_{m_1}) = - t_{m_1}^{-2} + 1 - t_{m_1}^2.\]

\subsection{Example 1}

\label{example1section}

We first produce an example of a nonsymplectic 4-manifold which admits a
free circle action whose orbit space is a
3-manifold which is fibered over the circle. Our construction generalizes
easily to produce a large class of such
manifolds with this  property.  Let $K_1$  and $K_2$ be any  fibered knots.
Form the fiber sum of the complements of
$K_1$ and $K_2$ with neighborhoods of the first and second meridians of
$T^3$, i.e.,
\[M_{K_1 K_2} = (S^3\setminus K_1)\#_{m=m_1} T^3 \#_{m_2=m} (S^3 \setminus
K_2) \]
where $m$ is the meridian of the corresponding knot.   Since both $K_1$ and
$K_2$ are fibered, the manifold
$M_{K_1 K_2}$ is a fibered 3-manifold. By Meng-Taubes theorem, the
Seiberg-Witten polynomial of
this manifold is
\[\mathcal{SW}^3_{M_{K_1 K_2}}(t_{m_1},t_{m_2}) =
\Delta_{K_1}(t_{m_1}^2)\Delta_{K_2}(t_{m_2}^2).\]
Let $X_{K_1 K_2}(l)$  be the 4-manifold
with free circle action that has $M_{K_1 K_2}$ for its orbit space and
$PD[l]$ for the Euler class of the circle
action. Taking both $K_1$ and $K_2$ to be the figure eight knot (see
Figure~\ref{fig:borromean2}), we get a
manifold with Seiberg-Witten polynomial:

\begin{figure}
\includegraphics{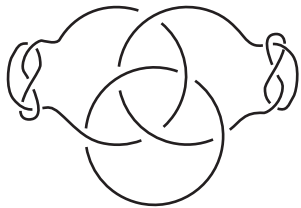}
\caption{$M_{K_1 K_2}$ before surgery }
\label{fig:borromean2}
\end{figure}

\begin{multline*} \mathcal{SW}^3_{M_{K_1 K_2}} = t_{m_1}^{-2}t_{m_2}^{-2} -
3t_{m_2}^{-2} + t_{m_1}^2 t_{m_2}^{-2} -3
t_{m_1}^{-2} + 9 \\ -3t_{m_1}^2 + t_{m_1}^{-2}t_{m_2}^2 - 3t_{m_2}^2 +t_{m_1}^2
t_{m_2}^2.\end{multline*}
The Seiberg-Witten polynomial of the manifold $X_{K_1 K_2}(4m_1)$ can be
calculated from
Theorem~\ref{thm:maintheorem},
\[\mathcal{SW}^4_{X_{K_1 K_2}(4m_1)} = 2t_{m_1+m_2}^{-2} -3 t_{m_2}^{-2}  +
9 - 6 t_{m_1}^2 +2t_{m_1 + m_2}^2
-3t_{m_2}^2,\]
where $t_\beta = \mbox{exp}({\pi^*(PD(\beta))})$ is the pullback of the
\Spinc\ structure on $M_{K_1 K_2}$.

A theorem of Taubes \cite{sw:sw_inv_and_symp_form} implies that the first
Chern class $c_1$ of a symplectic 4-manifold
must have Seiberg-Witten invariant $\pm1$. We thus see that the manifold
$X_{K_1 K_2}(4m_1)$ admits no symplectic
structure with either orientation. This is not the only free
$S^1$-manifold over $M_{K_1 K_2}$ with this property.  The manifolds
$X_{K_1 K_2}(-4m_1)$, $X_{K_1 K_2}(4m_2)$,
and $X_{K_1 K_2}(-4m_2)$ also admit no symplectic structures.

\subsection{Example 2}
Next we produce an example of a $3$-manifold which is not the orbit space
of any symplectic $4$-manifold with a free
circle action. Let $K_1=K_2$ be the nonfibered knot $5_2$ (see
~\cite{knots:knotsandlinks}).  The Seiberg-Witten
polynomial of $M_{K_1 K_2}$ is
\begin{multline*}
\mathcal{SW}_{M_{K_1 K_2}}^3 = 4t_{m_1}^{-2}t_{m_2}^{-2} - 6t_{m_2}^{-2}
+4t_{m_1}^2 t_{m_2}^{-2} - 6
t_{m_1}^{-2} + 9 \\ -6t_{m_1}^2 +4t_{m_1}^{-2}t_{m_2}^2 - 6t_{m_2}^2 +
4t_{m_1}^2 t_{m_2}^2.
\end{multline*}
One then needs to calculate as in Example 1. There are only finitely many
free $S^1$ manifolds $X_{K_1 K_2}(l)$ which
need to be checked because for all $l=a m_1 + b m_2$ with $|a|,|b| > 2$ the
Seiberg-Witten polynomial
$\mathcal{SW}^4$ is equal to the  3-dimensional polynomial  (only the
meaning of the
variables will change).  A calculation shows that the remaining free
$S^1$-manifolds all have \Spinc\ structures with  Seiberg-Witten invariant
greater than one in absolute value.
Therefore these manifolds are not symplectic.
Therefore $M_{K_1 K_2}$ is not the orbit space of  any symplectic
$4$-manifold with a free circle action.

\subsection{Remarks}

The above two examples show:

\smallskip
\noindent {\bf{1}}. {\it{There exist nonsymplectic free $S^1$-manifolds
with fibered orbit space.}}

\smallskip
\noindent {\bf{2}}. {\it
{There exists a $3$-manifold which is not the orbit space of any symplectic
4-manifold with a free $S^1$-action.}}

\smallskip
We conclude with two questions.

\begin{quest}
If $X$ is a free $S^1$-manifold which is symplectic, must its orbit space
$M = X/S^1$ be fibered?
\end{quest}
Taubes has conjectured this in case $X=M\times S^1$.
Theorem~\ref{thm:maintheorem} could be used to search for manifolds with
free $S^1$-actions that had nonfibered orbit
spaces and which do not have Seiberg-Witten obstructions to having
symplectic structures.  One would still need to prove
that  those manifolds where symplectic.  While a counter example may be
obtainable, a proof to the affirmative is
already at least as difficult as a proof of Taubes' conjecture.

\begin{quest} Let $M$ be a 3-manifold with the property that every free $S^1$-manifold whose orbit space is $M$ 
is
symplectic. Is $M$ fibered?
\end{quest}

The $3$-torus is an example of manifold with this property~\cite{symp:compact_parallel_symp_complex_man}.







\end{document}